\documentclass[onecolumn,notitlepage,aps,pra,10pt]{revtex4-1}
\usepackage{amsmath}
\usepackage{color}
\usepackage{graphicx} 

\begin{document}

\title{Ropelength-minimizing concentric helices and non-alternating torus knots}
\author{Alexander R. Klotz} 
\affiliation{Department of Physics and Astronomy\\ California State 
University, Long Beach \\ alex.klotz@csulb.edu}
\author{Finn Thompson}
\affiliation{School of Mathematics and Physics\\ University of Queensland \\ f.thompson@uq.net.au}

\begin{abstract}
    An alternating torus knot or link may be constructed from a repeating double helix after connecting its two ends. A structure with additional helices may be closed to form a non-alternating torus knot or link. Previous work has optimized the dimensions and pitch of double helices to derive upper bounds on the ropelength of alternating torus knots, but non-alternating knots have not been studied extensively and are known to be tighter. Here, we examine concentric helices as units of non-alternating torus knots and discuss considerations for minimizing their contour length. By optimizing both the geometry and combinatorics of the helices, we find efficient configurations for systems with between 3 and 39 helices. Using insights from those cases, we develop an efficient construction for larger systems and show that concentric helices distributed between many shells have an optimized ropelength of approximately $7.83Q^{3/2}$ where $Q$ is the total number of helices or the minor index of the torus knot, and the prefactor is exact and a 75\% reduction from previous work. Links formed by extending these helices and bending them into a $T(3Q,Q)$ torus link have a ropelength that is approximately 12 times the three-quarter power of the crossing number. These results reduce the ratio between the upper and lower bounds of the ropelength of non-alternating torus knots from 29 to between 1.4 and 3.8.
\end{abstract}
\maketitle

\section{Introduction}

The ropelength of a knot is the minimum contour length that a knot can have while respecting the no-overlap constraint of a tube of unit radius with the knot along its axis \cite{cantarella2002minimum}, or more simply put the least amount of rope required to tie a given knot. Ropelength connects the mathematics of knot theory with physical knots, and is known to be related to deeper topological invariants \cite{diao2020braid}. Besides the unknot and certain Hopf chain links, no knot or link has a known ropelength. There are generally three approaches to ropelength research. One is to rigorously establish upper and lower bounds on how ropelength scales with crossing number by relating relating it to other topological invariants, establishing geometric constraints on closed curves, or mapping the knots to lattice walks \cite{buck1999thickness, diao2003lower, diao1998complexity}. Another is to use optimization algorithms such as SONO \cite{pieranski1998search} or Ridgerunner \cite{ashton2011knot} to find tight upper bounds on the ropelengths of specific knots, examining the overall scaling of ropelength with crossing number or the properties of the tight knots themselves \cite{klotz2024ropelength}. Finally, certain knots or classes of knots can be constructed and optimized to establish upper bounds for the ropelength-crossing number relationship of these families, including petal knots \cite{tinlin}, 2-bridge knots \cite{huh2021tight}, and arborescent links \cite{mullins}. 

The construction method has been extensively applied to alternating $T(P,2)$ torus knots or links, for which a single unit can be defined and concatenated, with each additional unit adding a fixed contour length and two crossings. Bohr and Olsen treated the knots as circular double helices and arrived at the optimal helical pitch and radius, finding an asymptotic ropelength upper bound of 17 per twist or 8.5 per crossing \cite{olsen2013principle}. Huh, Oh, and Kim treated the knots as superhelices of four winding curves, finding 7.63 per crossing \cite{huh2018ropelength}. Klotz and Maldonado treated alternating torus links as circles wrapped by single helices, reducing the coefficient to 7.36 \cite{klotz2021ropelength}. Most recently, Kim et al. found the ideal radii of a larger helix wrapped around a smaller, finding a current-best value of 7.32 \cite{kim2024efficiency}. Although numerical evidence indicates that this is close to the true minimal value \cite{pieranski1998search}, the best proven lower bound is 1/28 per crossing \cite{diao2020braid} for links and 1/56 for knots \cite{diao2022ropelength}. 

Given the considerable effort devoted to bounding the alternating torus knot ropelength, we believe non-alternating knots deserve study as well. Non-alternating knots are known to be tighter than alternating knots, both from numerical optimizations \cite{klotz2024ropelength} and a proven sub-linear scaling with crossing number of at least a 3/4 power \cite{buck1999thickness}. Diao and Ernst \cite{diao1998complexity} showed that lattice constructions of non-alternating torus knots could achieve the 3/4 scaling, and Cantarella et al. used helical geometry to establish the same exponent \cite{cantarella1998tight} without explicitly constructing helical knots. For every crossing number at which there is a torus knot and all ropelengths have been measured (3, 5, 7, 8n, 9a, 10n, 11a), the tightest knot in each class is a torus knot, suggesting that torus knots are an ideal system for studying the applicability of ropelength bounds. Given that Diao and Ernst \cite{diao1998complexity} and Cantarella et al. \cite{cantarella1998tight} both demonstrated that the ropelength of a non-alternating torus knot is at least $L=\kappa C^{3/4}$, we would like to constrain what that constant $\kappa$ is. Based on methods and results from Bohr and Olsen \cite{olsen2012geometry, olsen2013principle} and from Kim, Oh, and Huh \cite{huh2016best}, we extend constructions of optimal helices to non-alternating torus links based on concentric helices. We do this by constructing ``multihelices,'' concentric shells of helices around a central rod, and optimizing the contour length based on the height and radii of the helices and the arrangement of helices in each shell. The number of helices plus the rod gives the smaller index, $Q$, of a torus knot or link. Since increasing $Q$ will increase both the ropelength and the number of crossings, the figure of merit is the ropelength per crossing. After optimizing individual helical units, we will discuss their closure into full torus knots and links. Ultimately, we will extend the values of $Q$ for which torus knots have been optimized from 2 to 39, and establish an upper bound for non-alternating torus links that scales with the same 3/4 power as the lower bound but is stronger than the upper bound established by Diao and Ernst using lattice methods \cite{diao1998complexity} .


\section{Cylindrical $N$-helices}

If $N$ helices with the same pitch and radius lie equally spaced around a circle in a plane perpendicular to the pitch of the helices, they can be said to wind around the exterior of a cylinder of height $H$ and radius $R$. An $N$-helix can be constructed from vertical lines arranged in a circle by applying a rotation matrix to each of the lines with an angle $\theta$ that increases from 0 at the bottom to $2\pi$ at the top. The equation of the $i$th helix in an $N$-helix is:

\begin{align}
    x(\theta)=R\cos\left(\theta+i\frac{2\pi}{N}\right),\ \ y(\theta)=R\sin\left(\theta+i\frac{2\pi}{N}\right),\ \ z(\theta)=\frac{H\theta}{2\pi}
\end{align}

Along every 2$\pi$ pitch of the helices, each helix passes over and under every other helix. The total number of crossings is $N(N-1)$, the factor of two from the crossings canceling the factor of two from double-count avoidance. The total arc length of each helix is simply the hypotenuse of the unwrapped cylinder:
\begin{equation}
    L_h=\sqrt{H^2 +(2\pi R)^2}
\end{equation}
For a helix representing a tube of unit radius, there is a no-overlap constraint on $R$, $H$, and $N$ such that no two helices can be within distance 2 of each other. Huh et al. derived a transcendental inequality based on the minimum distance between two neighboring helices \cite{huh2016best}. They found the conditions that minimize length, finding the ideal radius $R_o$ and height $H_o$. The transcendental equation that minimizes arc length is: 
\begin{equation}
    2-2\cos{\left(\phi_N+\frac{2\pi}{N}\right)}=\phi_N^2.
\end{equation}
The meaningful value of $\phi_N$ is the real negative solution to the above equation in the range $-\frac{2\pi}{N}i\leq\phi_N\leq0$. This represents the separation in phase between two helices at their closest approach. Finding this transcendental value allows the optimal radius, height, pitch angle, and ropelength to be determined based on the following expressions, values from which are shown in Table 1.
\begin{equation}
    R_o=2\sqrt\frac{1}{\phi_N^2-\phi_N\sin(\phi_N+\frac{2\pi}{N})},\ \ \ H_o=4\pi\sqrt\frac{\sin(\phi_N+\frac{2\pi}{N})}{\phi_N^2\sin(\phi_N+\frac{2\pi}{N})-\phi_N^3}.
\end{equation}
This radius and height, which depend only on $N$, are referred to as the ``ideal'' values even in conditions in which they do not globally minimize ropelength when multiple concentric helices are present. As $N$ increases, the terms within the cosine of the constraint equation decreases, such that it can be approximated by a Taylor series:

\begin{equation}
    \phi_N^2\approx2-2\left(1-\frac{1}{2}\left(\phi_N+\frac{2\pi}{N}\right)^2\right)=\phi_N^2+\frac{4\pi\phi_N}{N}+\frac{4\pi^2}{N^2 }
\end{equation}
This is solved by $\phi_N=\frac{\pi}{N}$, which can be substituted into the optimality equations to find simplified versions of the height, radius, ropelength, and pitch angle:

\begin{table}[ht]
\centering

\begin{tabular}{|l|l|l|l|l|}
\hline
$N$ & $R$ & $H$ & $L$ & $L/C$ \\ \hline
2 & 1.04587 & 5.3934 & 17.0026 & 8.5013 \\ \hline
3 & 1.43524 & 8.36102 & 36.8926 & 6.14877 \\ \hline
4 & 1.86156 & \textbf{11.2303} & 64.8603 & 5.40503 \\ \hline
5 & 2.29861 & \textbf{14.0787} & \textbf{100.846} & 5.0423 \\ \hline
6 & 2.7404 & \textbf{16.9191} & \textbf{144.839} & 4.82797 \\ \hline
7 & 3.18472 & \textbf{19.7556} & \textbf{196.834} & 4.68652 \\ \hline
8 & \textbf{3.63056} & \textbf{22.5898} & \textbf{256.832} & 4.58629 \\ \hline
9 & \textbf{4.07739} & \textbf{25.4227} & \textbf{324.83} & 4.51153 \\ \hline
10 & \textbf{4.52491} & \textbf{28.2546} & \textbf{400.828} & 4.45364 \\ \hline
\end{tabular}
\caption{Optimal radius ($R$) and height ($H$) for $N$ helices wound around a cylinder that minimizes the total length ($L$) of the helices. $L/C$ is the length per crossing of one helical unit. Bold values are described by the large-$N$ approximation within 1\%.}
\end{table}

\begin{align} 
R_o=\frac{\sqrt{2}}{\pi}N,\ \ H_o=\sqrt{8}N,\ \ \
L_o=4N^2
\end{align}

Although these approximations nominally apply in the large-$N$ limit, they are quite accurate for $N$ as small as 4, where the approximation overshoots the height by 0.7\% and undershoots the radius by 4\%. By $N$=12 the radius undershoot is below 0.4\%. We will use these approximations to derive expressions for the limiting cases of multihelical constructions, but when defining coordinates for specific helices or knots, the exact values should be used. Common features emerge independent of N: the circumference is always equal to the height and the pitch angle is always 45 degrees. In the large-$N$ limit where the crossing number is approximately $N^2$, it can be seen that an $N$-helix torus knot will have a ropelength of 4 per crossing, slightly more than half that of the best bound for alternating torus knots. This implies that the ropelength-per-crossing of a cylindrical N-helix cannot be improved beyond a certain point by adding more helices.

\begin{figure}[ht]
    \centering
    \includegraphics[width=0.8\linewidth]{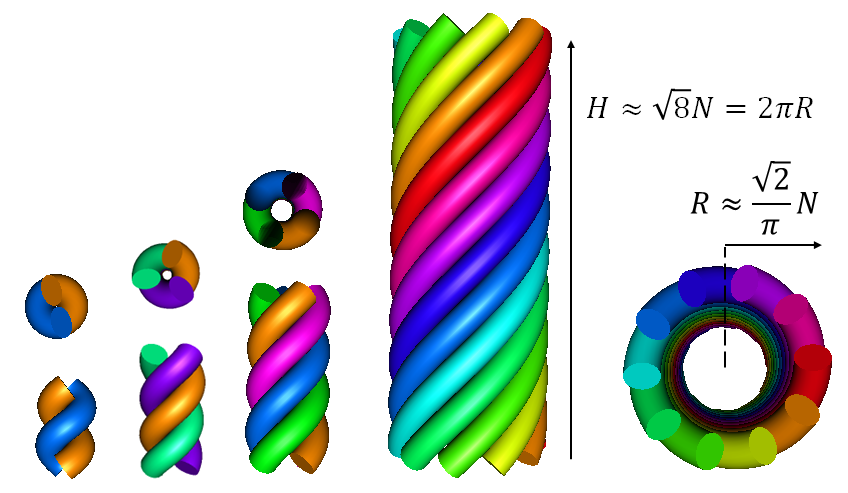}
    \caption{Side and top views of ideal cylindrical helices with 2, 3, 4, and 10 helices. The height and radius of the 10-helix are well approximated by the large-$N$ approximations.}
    \label{fig:helix}
\end{figure}

For more complex constructions, the height and radius are not necessarily at their ideal value. For a given height and number of helices, there is a minimum radius that can accommodate them, and for a given radius there is a minimum height. Suppose there exists tractable functions $\tilde{R}(N,H)$ and $\tilde{H}(N,R)$ that give the smallest height that can accommodate $N$ helices at a given radius, or the smallest radius at a given height. These functions would allow the optimization of more complex concentric helical structures. Unfortunately, there are no closed functional forms of these equations, but they can be described by a universal empirical function, presented in the Appendix.

\section{Asclepius and the Caduceus}

The Rod of Asclepius features a single snake winding around a rod, while the Caduceus features two snakes (Fig. 2a). In addition to their powerful iconography, they also represent efficient torus knot configurations.
A cylindrical $N$-helix can be made more efficient by straightening one of the helices and placing it in the center, which allows the radius and height of remaining $N-1$ helices to be reduced, or equivalently by adding a straight rod in the center of an $N$-helix to create a more efficient $N+1$-helix. Figure 2b shows a 3-helix and a more efficient 1+2-helix, which is shorter but wider. When $N\geq6$, the ideal radius of the $N-1$-helix is greater than 2 and the central rod does not interfere, and it is straightforward to show that converting one of the helices to a rod
will reduce the overall length. For $N$=2, 3, 4, and 5, the radius of the $N-1$-helix must be increased to 2 so as not to interfere with the central rod, and it is not obvious that the overall length will decrease. In each of these cases, the length decreases (Fig. 2c). The $Q=2$ Rod of Asclepius was demonstrated by Klotz and Maldonado (with the central curve being a circle rather than a rod) \cite{klotz2021ropelength}, and was improved slightly by Kim et al. by giving the rod a slight helical twist with a radius of 0.14 \cite{kim2024efficiency}. The large-$N$ approximation gives
\begin{equation}
    L_{cad}=4(N-1)^2+\sqrt8(N-1)<4N^2,
    \label{eq:cad}
\end{equation}
and describes $N\geq6$ quite well. Because of the required radius increase, there is less improvement from $N$=2-4 relative to the cylindrical helix, and the improvement is maximized for $N$=5. There is an additional improvement that can be made by increasing the radius of the helix beyond the ideal value to allow the height of the rod to decrease. It is extremely marginal and decreasing for $N\geq6$. For $N$=6, increasing the radius of the 5-helix from 2.29 to 2.37 decreases the height from 14.08 to 13.62, reducing the length of the unit from 114.94 to 114.57, a 0.3\% improvement. This improvement, comparable to half the width of a data point in Fig. 2c, is smaller for larger $N$.

\begin{figure*}
    \centering
    \includegraphics[width=1\linewidth]{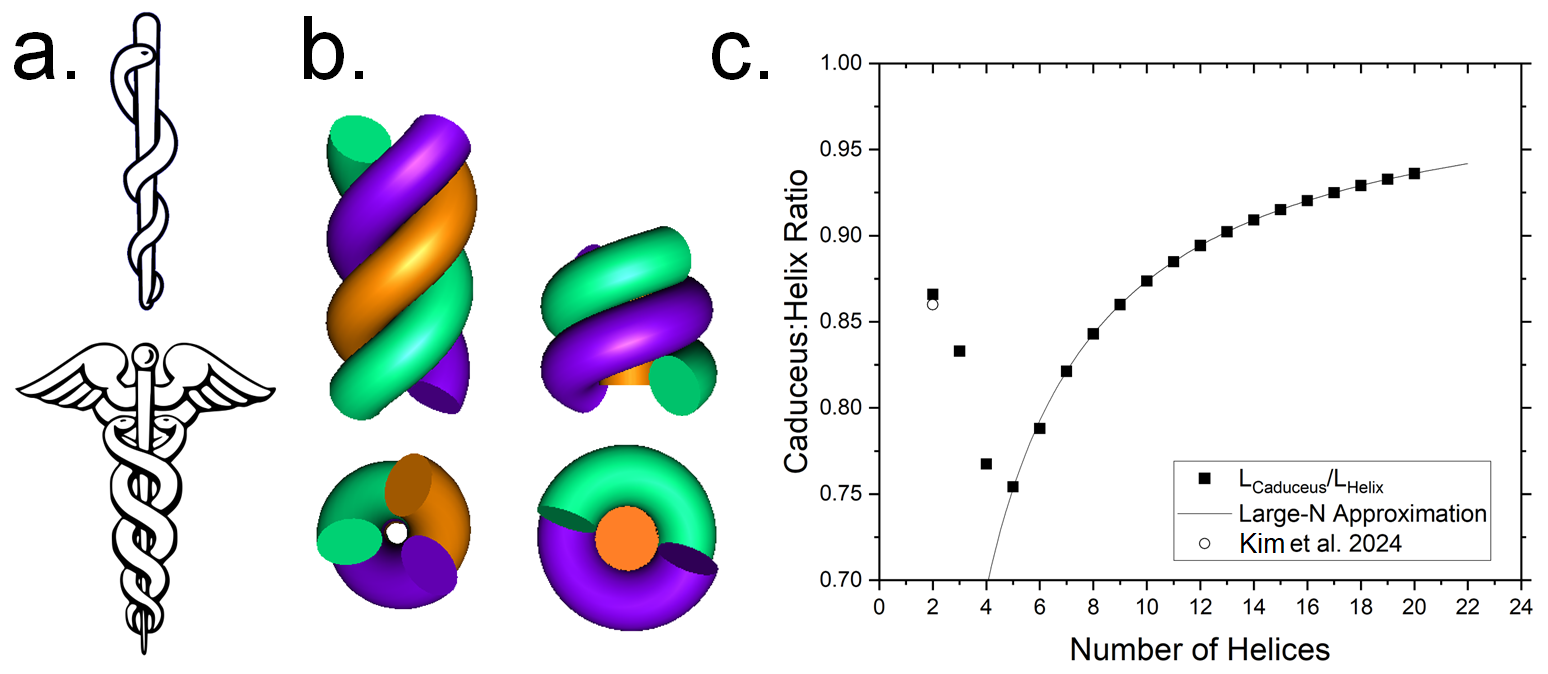}
    \caption{a. Artistic renderings of the Rod of Asclepius and the Caduceus, representing efficient 2- and 3-helices \cite{wiki}. b. An ideal triple helix, and a shorter but wider double helix wound around a central rod. c. Improvement of Caduceus-style helices relative to an $N$-helix with the same topology. Above 5 helices, the improvement is well-described by Eq. \ref{eq:cad}. For 2 through 5 helices, the radius must be increased, reducing the improvement.}
    \label{fig:cad}
\end{figure*}

\section{Optimizing concentric helical shells}
In the rest of this paper, we describe configurations of a central rod and multiple concentric rings of helices. We refer to these rings as shells, as the cross section of the helices bears resemblance to electron shell diagrams. A structure with a rod and $N$ helices (with $Q=N+1$) spread over $T$ total shells is referred to as a $Q$-multihelix. For a given number of helices there is some optimal arrangement of them into shells, but like electrons around a nucleus that arrangement can be non-trivial. Each shell with a given number of helices at a given radius has a minimum height $H$ at which the helices do not overlap. Since the top and bottom of the entire construct must be level, whichever of the $T$ total shells has the highest $H$ sets the height of the entire system, and the ropelength can be written as a sum over all the shells:
\begin{equation}
    L_{tot}=H+N_i\sum_{i=1}^T\sqrt{H^2+(2\pi R_i)^2 }.
\end{equation}
The challenge is finding the combination of $N_i$ and $R_i$  that minimize $L_{tot}$ for a given $Q$. The optimization is both combinatorial, determining the best distribution of helices between shells, and geometric, determining the best radii of the shells. The notation we use for the arrangement of a multihelix is 1-A-B-C..., e.g. the 1-4-5-2 multihelix depicted in Fig. 3 has a central rod, an inner shell of 4, a middle shell of 5, and an outer shell of 2 helices. 

We begin with a rod and two shells, of which $(N-M)$ helices are in the inner shell and $M$ are in the outer shell (Fig. 3). Although this is not the optimal configuration for $N>9$, it provides insights about optimizing ropelength that can be extended to more complex systems. Each shell must have sufficient radius such that no helix overlaps with another in its shell, and that the two shells are separated by a distance of at least 2. The height of the entire configuration is set by that of the most constrained helix. For a given $N$, there are three dimensions of optimization: the number of helices in the outermost shell $(M)$, the radius of the inner shell ($R_1$), and the radius of the outermost shell ($R_2$). 
It is useful to construct a configuration that does not require geometric optimization, height determination, or constraint checking. One such configuration (A) has the $(N-M)$ inner helices at their ``ideal'' radius $R_1=\frac{\sqrt{2}}{\pi}(N-M)$, as well as the M outer helices at $R_2=\frac{\sqrt{2}}{\pi}M$. The height of the system is set by the outer helix at $H=\sqrt{8}M$. The total length of this configuration is:
\begin{equation}
    L_A=\sqrt{8}\left(M+\sqrt{2}M^2+(N-M)\sqrt{M^2+(N-M)^2}\right)
    \label{eq:idealideal}
\end{equation}
Because the radii are constrained by $R_2-R_1\geq2$, the number of outer helices must be at least 5 greater than the number of inner helices, leaving a gap of approximately 2.25 between the two shells. The ropelength is minimized with respect to the one combinatorial degree of freedom when the outer shell has as few helices as possible, which is $M=(N-M)+5=(N+5)/2$. This yields an expression for the ropelength: 
\begin{equation}
    L_{A0}=N^2+\left(10+\sqrt{2}\right)N+(N-5)\sqrt{N^2+25}+25+\sqrt{50}
    \label{eq:ideal5}
\end{equation}
This expression has an asymptotic length of $2N^2$, or 2 per crossing,  which is to be expected if a single ideal $N$-helix were split into two shells with the same crossing number but half the length. More generally, for a 2-shell multihelix with a given $N$, we can minimize the total length with respect to the radii of each shell for each possible configuration. Fig. 3 shows four examples of a 2-shell 12-multihelix, for which the ropelength is minimized by the 1-6-5 configuration. We find, universally, that the distance between the two shells takes its smallest value of $R_2=R_1+2$, except in non-optimal arrangements (e.g. 2 inner helices and 9 outer helices). For every $Q$ from 6 to 26, the three most efficient configurations all have a radial separation of exactly 2 between the shells. The is also true for 3-shell configurations optimized for all three radii, for $Q$ from 6 to 26. Assuming this feature holds, it reduces the geometric dimensionality for a $T$-shell system from $T$ to 1.

Distributing helices between a larger number of shells often decreases the ropelength, but it is inefficient to have one shell for each helix. Finding the ideal arrangement for a large number of helices is nontrivial: the number of possible arrangements is the number of integer compositions of $N$, which is $2^{N-1}$.  By constructing a multihelix for each of the $2^{N-1}$ arrangements while varying the inner radius, we can exhaustively find the globally minimal configuration for each $N$. Beyond $N$=20 this becomes computationally challenging, and we only search among configurations with a number of shells likely to produce an optimal configuration, e.g. within $\pm$1 of the $(N-1)$th.

\begin{figure*}[ht]
    \centering
    \includegraphics[width=1\linewidth]{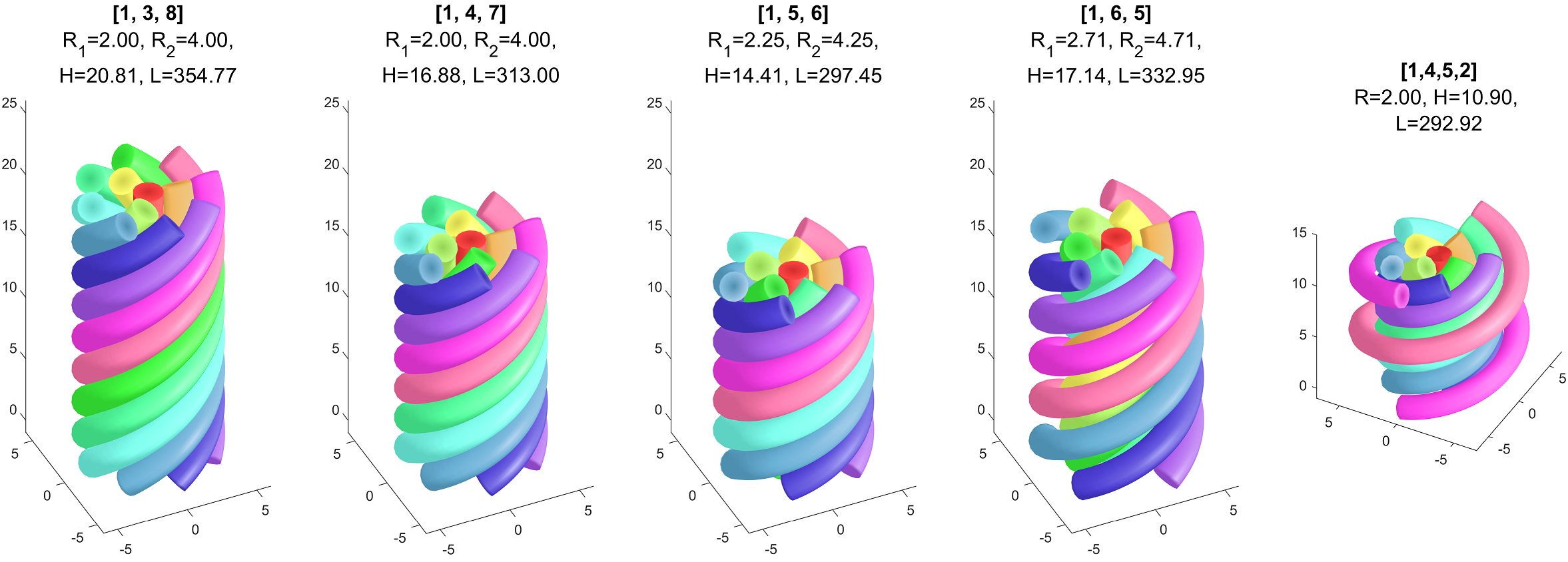}
    \caption{Four geometrically optimized 12-multihelices with a rod and a varying distribution of helices between the first and second shells. The contour length is minimized for the 1-5-6 configuration. In all cases, the second shell is radially separated from the first by a distance of 2. The fifth image is the geometrically \textit{and} combinatorially optimized 12-helix, which has a rod and three helices in a 1-4-5-2 arrangement at radii 2, 4, and 6. } 
        \label{fig:triple}
\end{figure*}
With each shell being separated from its neighbors by a distance of 2, the total length of a configuration of $Q=N+1$ helices over $T$ shells is
\begin{equation}
L=H+\sum_{i=1}^T M_i \sqrt{H^2+(2\pi (R+2(i-1)))^2}\end{equation}
where $M_i$ is the number of helices in the $i$th shell, and $\sum_{i=1}^T M_i=N$. Here, $H,R$ are the optimal values that minimize $L$ while satisfying the no-overlap constraint. 
With a central rod, the first shell is at a radius $R\geq2$, and all other shells are tightly packed at radii $R+2,R+4,\cdots$. To avoid overlaps with the central rod and between the cross-sections of the helices, we must have $R\geq\max\{2,1/\sin(\pi/M_1)\}$. Without loss of generality, we may assume $M_1\neq0$. Similarly, a tight configuration will not have any empty inner shells, which requires that $R\leq4$. Furthermore, we require $H\geq 2M_i$ for each $M_i$.

This is ultimately a nonlinear optimization problem. We note that it is possible to solve for all of $H,R,T,M_i$ for a given $N$, but this involves solving a mixed integer nonlinear system of up to $3+T$ variables. To ensure optimality, we instead compute ropelength minimizing values of $H,R$ for prescribed choices of $M_i$ for each $N$. In particular, we form the following system, 
\[2\leq R\leq4,\quad H\geq2\max\{{M_i}\},\quad \min\{\mathrm{dist}(h(R_i,H,0),h(R_i,H,\frac{2\pi}{M_i}))\}\geq2\,\,\forall i,\]
where $h(R_i,H,\theta)$ represents a helix of radius $R_i$ and height $H$ with an angular phase shift of $\theta$. The minimum distance is computed by generating discrete helices with 500 vertices for two adjacent helices in each shell, and calculating the minimum of the distance between all pairs of points. We note that if $M_i=1$ (only one helix) in a particular shell, the minimum distance constraint is automatically satisfied. With the distance function constraint and the bounds on $R,H$, we used MATLAB's \textit{fmincon} function, which uses the interior-point method \cite{byrd}, to find the best configuration for a given $N$ across all $2^{N-1}$ arrangements.

Some examples of ropelength minimizing values of $R$ and $H$ for various configurations of 12 helices are given in Fig. \ref{fig:triple}. The ropelength of 12 helices in 2 shells is minimized with the arrangement 1-5-6, having length $L=297.45$, with the first shell being placed at radius $R_1=2.25$, and the second shell at $R=4.25$. Exploring configurations with more than 2 shells yields the geometrically and combinatorially optimized 12-helix with the arrangement 1-4-5-2, with a length of $L=292.92$.

Some of the results of these optimizations are summarized in Table II, with a complete listing in the Appendix. The optimized ropelength per crossing can be seen in Fig. \ref{fig:bigboi} and empirically decreases roughly as $7.6/\sqrt{Q}$. Generally speaking, helices are added onto the outermost shell until it becomes more efficient to add a new shell. The ideal number of shells increases from 1 to 2 at $N$=4, 2 to 3 at $N$=10, 3 to 4 at $N$=21, and 4 to 5 first at 29, going back to 4 for $N=30$ and 31, and back to 5 at $N$=32. In most cases the height of the system is set by the second or third shell, which in larger systems allows many shells that are ``too wide'' at larger radius rather than ``too tall'' at smaller radius. It is inefficient to have shells with very different height constraints, as the ropelength could be reduced by moving helices from a constrained to an unconstrained shell. In most cases, the radii of the shells shell are simply 2, 4, 6, 8...meaning that only combinatorial and not geometric optimization is required (assuming that the height of each shell is known). In a few cases, the inner radius was 2.036. These cases were arrangements of the form 1-5-7... that were set by the second shell.

It is a reasonable hypothesis that a dense circular packing of circles, when extended and twisted, would produce an optimal helical construct, similar to analysis of ideally packed helical bundles \cite{atkinson2019constant}. However, such configurations are inefficient because the inner shells are underconstrained and the construct can be improved by moving helices from the outer shells to the inner. The 1-6-12 configuration is at least the 2198th most efficient, and 1-6-12-18 is not in the top 30,000. Helices need not be arranged uniformly in concentric circles, and Cantarella et al. \cite{cantarella1998tight} used 16 helices in a 4x4 square pattern as a demonstration. Although it was not their intent to present this as an optimal configuration, we found that such a construction is minimized when the helices are each separated from their neighbors by a distance of 2.85 and have a height of 26.75, giving a ropelength of 619.4. This is 35\% greater than our best concentric 16-multihelix.

\begin{table}[ht]
\begin{tabular}{|l|l|l|l|}
\hline
\begin{tabular}[c]{@{}l@{}}Number of\\ Components\end{tabular} & \begin{tabular}[c]{@{}l@{}}Ropelength\\ per Twist\end{tabular} & \begin{tabular}[c]{@{}l@{}}Ropelength\\ per Crossing\end{tabular} & \begin{tabular}[c]{@{}l@{}}Ideal\\ Arrangement\end{tabular} \\ \hline

1+5 & 101.8831 & 3.918581 & {[}1,\textbf{3},2{]} \\ \hline

1+10 & 253.6763 & 2.511646 & {[}1,4,\textbf{5},1{]} \\ \hline

1+15 & 459.9192 & 2.035041 & {[}1,5,\textbf{7},3{]} \\ \hline

1+20 & 674.2664 & 1.681462 & {[}1,5,7,\textbf{8}{]} \\ \hline

1+25 & 940.6767 & 1.502678 & {[}1,5,7,\textbf{8},5{]} \\ \hline
1+30 & 1226.158 & 1.360886 & {[}1,5,\textbf{8},9,8{]} \\ \hline
1+35 & 1545.313 & 1.260451 & {[}1,5,\textbf{8},9,9,4{]} \\ \hline
\end{tabular}
\caption{Abridged table of ropelength (per twist and crossing) of optimal concentric helices, with the ideal arrangement of helices. The bold digit indicates which shell sets the height of the system. The 1+$N$ indicates a central rod surrounded by $N$ helices. Full table is found in the Appendix.}
\end{table}

\section{Highly Concentric Helices}
Although it is not feasible to exhaustively find the ideal arrangement of a large number of helices, we can construct efficient arrangements to look at the limiting ropelength trends as that number increases. The crossing number per helical twist is approximately $Q^2$, and we will examine how the ropelength per crossing varies with increasing $Q$. We will derive increasingly efficient bounds on the length of a multihelix as a function of $Q$, and all will scale as $Q^{3/2}$ which is equivalent to the ropelength per crossing scaling as $Q^{-1/2}$. The ropelengths per crossing we derive hold for $T(P,Q)$ knots and links with large $P>Q$, and the effect of closure will be discussed in the next section.

\begin{figure*}
    \centering
    \includegraphics[width=1\linewidth]{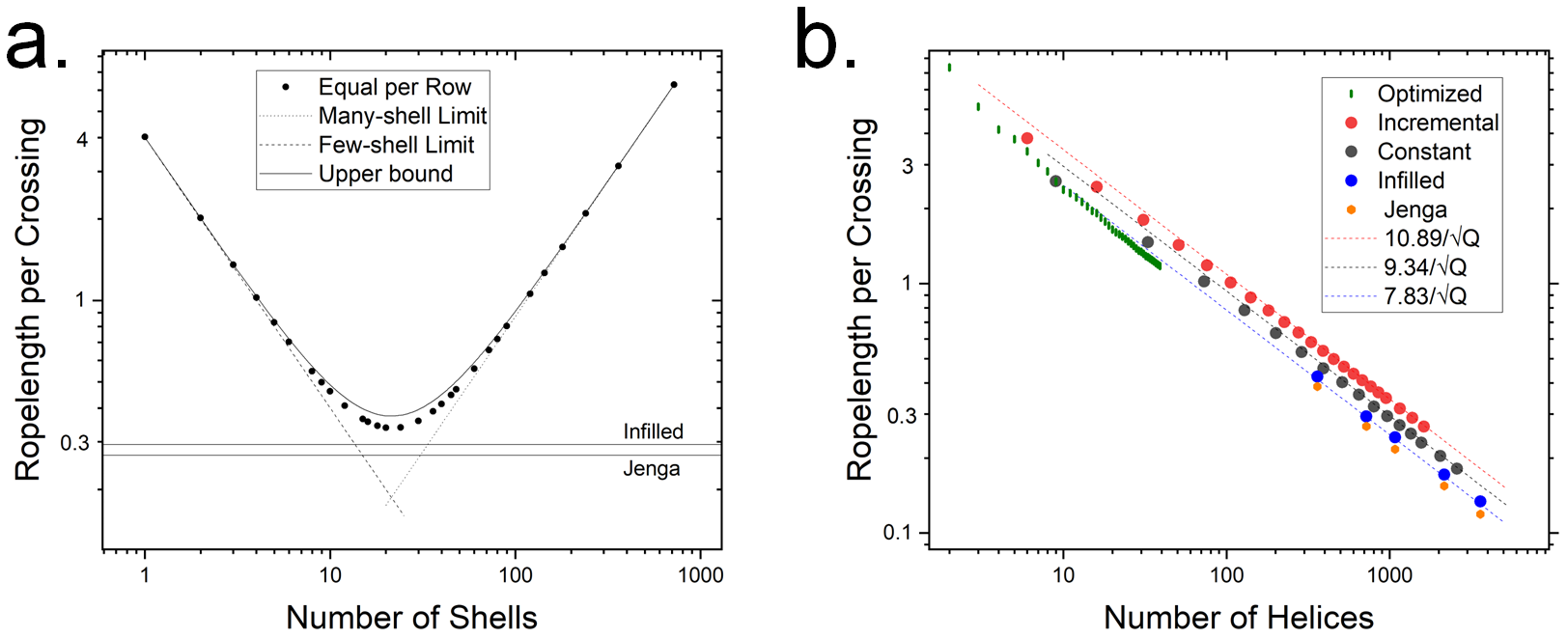}
    \caption{a. Ropelength per crossing of a winding of a 721-multihelix with a central rod and an equal number of helices per shell, as a function of the number of shells. The dependence is well-described by an upper bound derived from the sum of two limiting cases, Eq. \ref{eq:upper}. The horizontal lines show the improvements due to infilling the interior of the system, and further optimizing the arrangements by reverse-Jenga. b. Ropelength per crossing as a function of the total number of helices (including the rod). Green ticks show optimized configurations determined by exhaustive combinatirial sweeps and geometric minimization. Red data corresponds to an increment of 5 helices per shell, and the limiting trend from Eq. \ref{eq:1089}. Black data corresponds to an equal number of helices per shell, and the limiting trend from Eq. \ref{eq:934}. Blue data corresponds to infilling the equal-per-shell helices in increments of 4, and the limiting trend from Eq. \ref{eq:784}. Orange points show further improvements from a reverse-Jenga algorithm. All limiting trends depend inversely on the square root of the number of helices.}
    \label{fig:bigboi}
\end{figure*}

A construction that does not require optimization or constraint-checking involves incrementing the number of helices in each shell by a constant $k\geq5$, with each shell placed at its ideal radius, and the outer shell dictating the height of the system. For example, a 1-5-10-15-20 multihelix would have each shell at radius 2.29, 4.50, 6.75, and 9.00, a total height of 56.57, and a ropelength per twist of 3001.7. The ropelength of such a construction with T total shells may be written as:
\begin{equation}
L_{\alpha}=\sqrt{8}kT+4k^{2}T^{2}+\sqrt{8}k^2\sum_{i=1}^{T-1}i\sqrt{T^{2}+i^{2}}\approx1.72k^2T^3.
\label{eq:alpha}
\end{equation}
For large $T$, the dominant term is the sum, and its convergent value is derived in the Appendix. The total number of helices $Q$ is $\frac{k}{2}T(T-1)+1\approx\frac{k}{2}T^2$ and $T=\sqrt{2Q/k}$, meaning the ropelength per crossing is asymptotically:
\begin{equation}
    \frac{L_{\alpha}}{C}=\frac{1.72k^2T^3}{(\frac{k}{2}T^2)^2}=\frac{6.88}{T}=6.88\sqrt{\frac{2}{kQ}}\geq\frac{10.89}{\sqrt{Q}}
    \label{eq:1089}
    \end{equation}
The inequality is bound by the minimum value of $k$, 5, which leaves a gap of at least 2 between adjacent shells. This shows that the ropelength per crossing continually decreases as more shells are added, and can be closed into a link with the same asymptotic scaling with crossing number as the proven lower bound ($C^{3/4}$). However, the efficiency of large helices can be improved significantly. Following insights from the two-shell helices and the global minimization, we can surmise that an ideal arrangement likely has a similar number of helices per shell ($N_s$), that each shell except the innermost is the minimum distance from its neighbors, and one of the innermost shells sets the height of the system. With these considerations we can construct a multihelix with an equal number of helices per shell, with one combinatorial and one geometric degree of optimization. For a highly composite number of helices, 24 for example, a multihelix can be constructed with 1, 2, 3, 4, 6, 8, 12, or 24 shells, accordingly with 24, 12, 8, 6, 4, 3, 2, or 1 helix in each shell. The innermost shell is the most constrained and sets the height of the system, and each subsequent shell has radius 2 beyond its inner neighbor. The radius of the inner shell ($R_s$) can be optimized by computing the ropelength of the entire construct using $\tilde{H}$. An example for 720 helices about a central rod is shown in Figure \ref{fig:bigboi}a, where the ropelength is minimized at 20 shells of 36 helices each. 

The ropelength of such a construct may be written as:
\begin{equation}
    L_{\beta}=\tilde{H}(N_s,R_s)+N_s\sum_{i=1}^{T}\sqrt{\tilde{H}(N_s,R_s)^2+\left(2\pi\left(R_{s}+2(i-1)\right)\right)^2}
\end{equation}
We can estimate this dependence by examining the limiting cases. When the total number of shells $T$ is large, the radius of the outer helices far exceeds their height as well as the innermost radius, and only the $i$-dependent terms in the radical contribute:
\begin{equation}
    \frac{L_{\beta T}}{C}\approx\frac{1}{Q^2}\left(N_s\sum_{i=1}^T\sqrt{(2\pi( 2(i-1)))^2}\right)\approx \frac{1}{Q^2}\left(4\pi N_s\frac{T(T+1)}{2}\right) \approx\frac{2\pi T}{Q}
\end{equation}
When $Q$ is large but the number of shells is small, the radius of the inner shell far exceeds the distance of 2 between each subsequent shell. We can take the ``ideal'' height and radius of a single shell, where $R=H/2\pi$, to construct the $N_s$-dominant expression:

\begin{equation}
    \frac{L_{\beta N_s}}{C}\approx\frac{1}{Q^2}\left(N_s\sum_{i=1}^T\sqrt{8N_s^2 +(2\pi\sqrt{8} N_s/2\pi)^2 }\right)=\frac{1}{Q^2}\left(N_s\sum_{i=1}^T\sqrt{16N_s^2}\right)=\frac{1}{Q^2}\left(4TN_s^2\right) =\frac{4}{T}
\end{equation}
This is consistent with our earlier results that a single shell will reach a ropelength per crossing of 4, and that a double shell will reach 2. For intermediate values of $N_s$ and $T$, the ropelength will not exceed the sum of the two limiting cases (see Appendix), allowing an upper bound to be derived:

\begin{equation}
    \frac{L_{\beta}}{C}<\frac{2\pi T}{Q}+\frac{4}{T}
    \label{eq:upper}
\end{equation}
This upper bound, overlaid in Fig. \ref{fig:bigboi} describes the general trend of the data quite well. The upper bound is minimized at $T=\sqrt{\frac{2Q}{\pi}}$ which leads to an optimized expression for the ropelength of this construction:
\begin{equation}
    \frac{L_{\beta0}}{C}<4\sqrt{\frac{2\pi}{Q}}\approx\frac{10.03}{\sqrt{Q}}
\end{equation}

Although the data in Figure \ref{fig:bigboi}a is based on optimizing the height and radius of the inner shell, if optimization is skipped and the ideal helix values are used, the helical construct is typically only 1\% longer. This upper bound can inform the derivation of a stronger expression. Based on the inner shell having the ideal radius and height and each subsequent shell being 2 greater in radius, and taking the optimal number of shells from the upper bound, we can write an expression for the length of this construction:

\begin{equation}
    L_\gamma=2\sqrt{\pi Q}+\sqrt{\frac{\pi Q}{2}}\sum_{i=1}^{\sqrt{2Q/\pi}}\sqrt{4\pi Q+\left(2\pi\left(\sqrt{\frac{Q}{\pi}}+2(i-1)\right)\right)^2}\approx9.34Q^{3/2}=\frac{9.34}{\sqrt{Q}}C
    \label{eq:934}
\end{equation}
The evaluation of this sum is found in the Appendix, and it should be understood that specific implementations of sums like these should use the closest integers to the radicals when appropriate. This is approximately 15\% more efficient than the incremental organization. Because these configurations have a significant amount of empty space between the rod and the first shell of helices, they can be further optimized by removing helices from the outermost shells and placing them inside the innermost shell, reducing their radius. If too many are placed in a single shell it risks increasing the height, but if each inner shell is underconstrained they can be filled  without changing the height.

If the innermost shell of the above construction has its ideal radius of $R_s=\sqrt\frac{Q}{\pi}$, we can add $m$ shells at $r_1=2$, $r_2=4$...$r_m=2m\leq R_s-2$. The total number of infilled shells is approximately $R_s/2$. If the $i$th shell contains $4i$ helices, they will be underconstrained as the ideal radius is $1.8i$, and the height is set by the shell at $R_s$. Assuming $Q$ is large, we can write an expression for the total number of helices that can be added to the interior:
\begin{equation}
    N_{int}=4+8+12 ... \frac{4}{2}\sqrt{\frac{Q}{\pi}}=4\left(1+2+3 ...\frac{1}{2}\sqrt{\frac{Q}{\pi}}\right)\approx\frac{Q}{2\pi}.
\end{equation}
These helices, roughly one-sixth of the total, are removed from the exterior and brought to the interior which greatly reduces their length. The number of removed shells is:
\begin{equation}
    \Delta T=-\frac{Q}{2\pi}\sqrt\frac{2}{Q\pi}=-\sqrt{\frac{Q}{2\pi^3}}
\end{equation}
We can compute the change in ropelength of the entire construct by subtracting the length of the removed helices and adding that of the new inner helices:
\begin{equation}
    \Delta L=-\sqrt\frac{\pi Q}{2}\sum_{i=\sqrt{2Q/\pi}-\sqrt{Q/2\pi^3}}^{\sqrt{2Q/\pi}}\sqrt{4\pi Q+\left(2\pi\left(\sqrt\frac{Q}{\pi}+2(i-1)\right)\right)^2}+\sum_{i=1}^{\frac{1}{2}\sqrt{Q/\pi}}4i\sqrt{4\pi Q+4\pi^2i^2}\approx-1.5Q^{3/2}
\end{equation}
These sums are evaluated for large $Q$ using the same integrals as are used for $L_\alpha$ and $L_\beta$ with the appropriate limits, and are -2.11$Q^{3/2}$ and +0.60$Q^{3/2}$ respectively. This gives an improved asymptotic expression for the length of a large helix:
\begin{equation}
    L_\delta=L_\gamma+\Delta L=7.83Q^{3/2}=\frac{7.83}{\sqrt{Q}}C.
    \label{eq:784}
\end{equation}
The prefactor of approximately 7.83 in an exact transcendental number that arises out of the integral solutions to the sums that determine the ropelength. It is:
\begin{equation}
\frac{L_\delta}{Q^{3/2}}=
\frac{1}{12\sqrt{\pi}}\left(
\begin{aligned}
    20\sqrt{5}-32-6\pi+\sqrt{1-(4+\sqrt{2})\pi+(5+\sqrt{8})\pi^2}\left(12\sqrt{2}-\frac{6}{\pi}\sqrt{2}+6\right) \\-3\sqrt{2}\pi\log{\left[\frac{2+\sqrt{2}}{\pi}\left(\sqrt{1-(4+\sqrt{2})\pi+(5+\sqrt{8})\pi^2}-\pi\left(1+\frac{1}{\sqrt{2}}\right)+1\right)\right]}  
\end{aligned}
\right)\approx 7.82869 
\end{equation}

We can compute exactly the correction for the 721-multihelix in Fig. \ref{fig:bigboi}, where the ideal configuration has 20 shells of 36, the innermost has a radius just over 16, allowing 7 shells with 112 helices to fit inside. This removes three full shells and an four additional helices from the outside, reducing the ropelength from 175,000 to 152,000. Further optimization is possible, but infilling increments of 5 instead of 4 runs into issues with the height constraint. While Eq. \ref{eq:784} is the best limiting expression in this work, helices can be optimized on a case-by-case basis through a reverse-Jenga algorithm. In the game of Jenga, blocks are removed from the bottom of a tower and placed on the top. In our reverse-Jenga algorithm, helices are removed from the outermost shell and placed in the innermost shell that they can be added to without inducing overlaps, which are checked by generating discrete helices. This proceeds until no more helices can be added to any shell. The initial configuration for such a procedure is the equal-per-shell after infilling in increments of 4 as described previously, but the algorithm works for any non-overlapping initial configuration. The innermost shell before infilling sets the height of the system according to its ``ideal'' value. The rod-and-720 helix system was found to have an optimal configuration of [1-5-12-17-22-26-30-33-36-38-40-41-42-43-44-44-45-46-46-47-47-16]. Infilling reduced the ropelength by 14\%, while reverse-Jenga optimization reduced it a further 9\%. We do not perform further geometric optimization, so better configurations likely exist. Instead of removing helices from the outermost shells, they can be removed from whichever shell is setting the height of the system, such that the overall height of the system can decrease after each Jenga move, or from all shells uniformly. We are unaware if this is more efficient.

Ropelength per crossing data from each of the discussed constructions is shown in Fig. \ref{fig:bigboi}b, including the globally minimal configurations determined by exhaustive combinatorics and geometric optimization in Section IV. Each construction follows its asymptotic trend beyond several hundred helices, corresponding to tens of thousands of crossings per twist. Notably, the reverse-Jenga constructions computed between 360 and 3600 helices appear to fall along the same trend as the globally optimized constructions below 39 helices, suggesting that the largest constructions are close to optimal.



\section{Closure Considerations}




There are two ways that repeating helical units are typically concatenated into closed knots or links. One is to wrap the entire configuration into a torus (Fig. 5ab), the other is to have two parallel helical constructs with connectors at either end. The curvature introduced by the torus method can violate the no-overlap constraint unless the major radius is increased slightly, while the connectors add a finite amount of ropelength and may or may not contribute to the crossing number. In the case where a large number of units are concatenated, the limiting ropelength behavior is not significantly affected by either closure mechanism.

A torus made from $p$ repeating units of the concentric helix closed upon itself yields a $T(pQ,Q)$ torus link, with $Q$ distinct components. We use lower case $p$ to distinguish from the capital parameter $P=pQ$ in a $T(P,Q)$ torus knot. As the helix of repeats $p$ times around the circumference of the circle, the major radius is $R_M=\frac{Hp}{2\pi}$. When the helical axis is bent into a circle, the inner extrema of the helices move closer together into an overlapping position. To restore them to a non-overlapping configuration, the major radius can be increased by the minor radius of the outmost helix, which puts their inner extrema back into a non-overlapping configuration. If the outer helix has the optimal radius $R_o=\sqrt{2}N/\pi=\frac{H}{2\pi}$, this is equivalent to increasing the major radius by a factor of $(p+1)/p$, or replacing $p$ with $(p+1)$ in the major radius.

\begin{figure}[ht]
    \centering
    \includegraphics[width=0.8\linewidth]{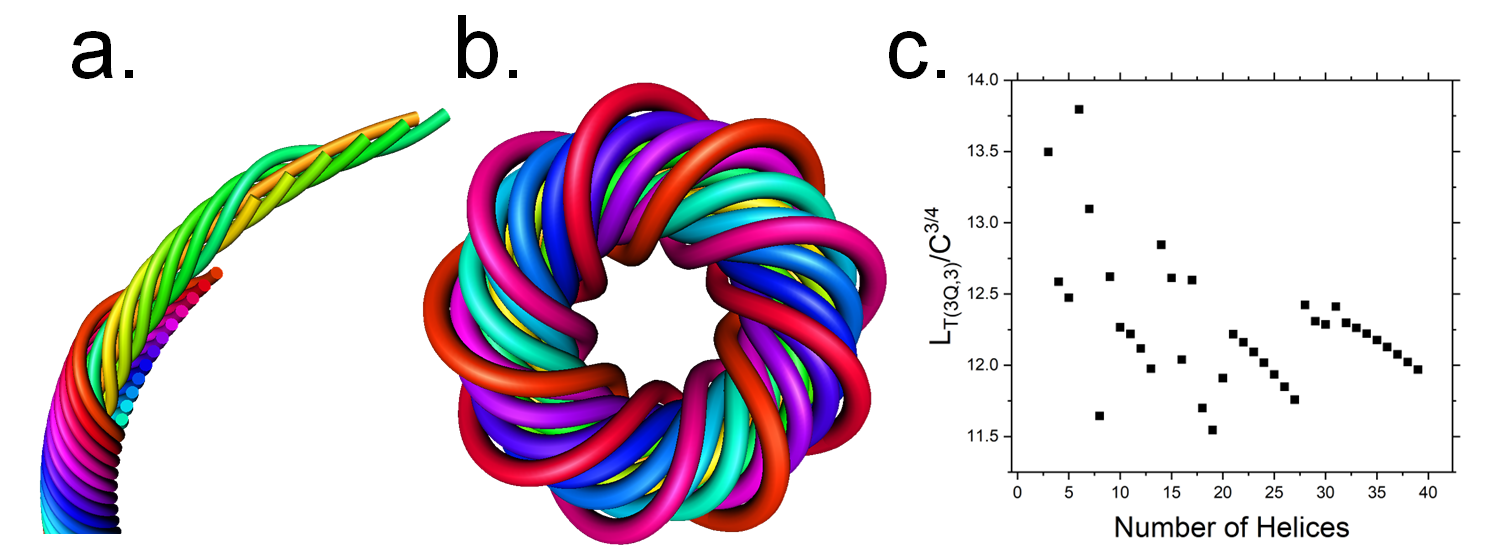}
    \caption{a. Partial rendering of a $T(266,19)$ torus link with a 1-6-12 configuration. b. Rendering of a $T(48,16)$ torus link with an optimized 1-5-7-3 configuration. c. Ratio of ropelength to the three-quarter power of crossing number of $T(3Q,Q)$ torus links constructed from optimized helical configurations.}
    \label{fig:closed}
\end{figure}

As the additional length due to closure proportionally decreases with increasing $p$, the ropelength per helix approaches the value of a single helical unit as $p$ increases. However, another figure of merit for non-alternating knots is $L/C^{3/4}$. This ratio is proven to be above 1.105 \cite{buck1998four}, and the smallest published value of a specific knot is 10.76 \cite{klotz2021ropelength}. Taking our incremental construction in which the outer helices are most constrained, assuming the increase in major radius required for closure scales the knot linearly, the crossing number is $pQ(Q-1)$, and the ratio is approximately:
\begin{equation}
    \frac{L}{C^{3/4}}\approx\left(\frac{p+1}{p}\right)\frac{10.89 pQ^{3/2}}{(pQ^2)^{3/4}}=10.89 p^{1/4}\left(\frac{p+1}{p}\right)
\end{equation}
This is minimized for $p=3$ and increases the numerical prefactor by a factor of 1.76 to 19.11. While this is an improvement upon lattice constructions \cite{diao1998complexity}, it is an overestimate for the ratio because the larger helices scale sublinearly with major radius. Since the outer helices in our more efficient constructs are underconstrained, the radius need not be increased as much for closure. Although $p=3$ was found to be optimal based on inefficient incremental constructions, $p=3$ was observed to be the best for other constructions as well. An optimized 1-5-7-6 helix can be closed into a $T(57,19)$ link with a ratio of 11.54 such that $2093=11.54\times1026^{3/4}$. Values of this ratio for optimized $T(3Q,Q)$ links are shown in Fig. \ref{fig:closed}c. The lower values in that figure are below the ratio found for any numerically tightened knot below 13 crossings, and are surpassed by numerical tightenings of $T(Q+1,Q)$ torus knots with at least 224 crossings \cite{klotz2021ropelength}. A conjectural lower bound can be derived from the fact that each of the $Q$ components has a linking number of 3 with the other $(Q-1)$ components, and thus cannot have a smaller arc length than the minimum convex hull around $3(Q-1)$ unit disks. In the large-$Q$ limit, a circular hull around a hexagonal close packing can be assumed, leading to the ratio $L/C^{3/4}=\sqrt{8\pi}\approx 5$, explained in further detail in the Appendix. Our overestimate on the upper bound of the ratio, 19.11, is a factor of 3.8 above this lower bound, a reduction from the factor of 29 between the upper bound for lattice torus knots of 32 \cite{diao1998complexity} and the lower bound of 1.1 \cite{buck1999thickness}. In practice, the convex hull around a smaller number of disks is not necessarily circular, and using the best known convex hulls \cite{cantarella2002minimum,kallrath2021near} we find lower bounds between 7 and 9.5 for the range of data shown in Fig. 5c. The lowest ratio of constructed upper bound to hull-based lower bound was 1.41, for the $T(93,31)$. Increasing the complexity does not necessarily decrease this ratio. The 1,557,360-crossing $T(2163,721)$ link constructed from the optimized 721-multihelix reaches a ratio of 12.06, which may be the most complex specific knot discussed in the literature to date. A previous analysis of $T(Q+1,Q)$ torus knot ropelengths using Ridgerunner \cite{klotz2021ropelength} showed that they scaled as the 0.71 power of crossing number. Despite being below the proven 0.75, this was interpreted as a ``finite size'' effect of the knots having not yet reached their asymptotic limit. The generally decreasing trend in Fig. 5c corroborates this, with a best-fit exponent to the ropelength of the $T(3Q,Q)$ links against crossing number of 0.73.

The previously reported upper bounds on the ropelength of non-alternating torus knots were from Diao and Ernst \cite{diao1998complexity}, who used a cubic lattice construction to find the most crossings possible for a lattice walk of a given length. This can be inverted into a constraint on a length required for a given crossing number. Converting to our notation, they found a 3/4-power bound with a coefficient between 20 and 32, not precisely defined. For any specific $Q$ we can create a configuration with a coefficient of approximately 12, and have established an upper bound of 19.11. Our results represent a 60 to 75\% reduction of the linear upper bound, and a 40 to 60\% reduction in the three-quarters bound.

Cartesian coordinates of these torus links can be constructed, although discretization may lead to spurious overlaps in tight configurations. The knot tightening software Ridgerunner \cite{ashton2011knot} detects overlaps in piecewise linear curves and rescales the coordinates accordingly, and can serve as an independent check for whether our constructions are non-overlapping. We generated a $T(48,16)$ link as three windings of a 1-5-10 concentric helix with 199 vertices per helix and input it to Ridgerunner. The initial contour length was rescaled from 2139.7 to 2141.6. With 100 vertices per helix, the rescaling was to 2146. We consider this evidence that the toroidal closures of our helices are non-overlapping, with the minor rescaling required due to discretization. Ridgerunner can also be used to assess the tightness of a configuration, although doing so is not the goal of this work. The aforementioned link reached a plateau at 1480.55, which is a reduction of 30\% compared to the non-optimized initial configuration. More efficient constructions are not necessarily better initial conditions for Ridgerunner. A more efficient 1-5-7-3 construction (Fig. 5b) with the same number of vertices reached a plateau at 1572. The algorithm would need to escape these local minima by inducing a buckling transition and breaking the symmetry of the link.

Creating a single-component knot from two parallel multihelices is more difficult. If each component on both constructs is labelled from 1 to $Q$, the $n$th component from the first construct can be connected to the $n$th at the same end of the second, and at the other end the $n$th component of the second can be connected to the $n+1$th of the first. This adds $(Q-1)$ essential crossings and results in a single component $T(pQ+Q-1, Q)$ torus knot. Since each of the $Q$ helices on multihelix must be connected at the top and bottom to the $Q$ helices on the other, and the smallest possible connection is a line segment of length 2, the additional length is at least $4Q$. Since the additional length also adds crossings, and the ropelength grows as $Q^{3/2}$, the closure length becomes negligible for large knots. However, if (for example) the connections requires a series of arcs whose required radius increases with $Q$, the scaling of the connector's length may not be negligible. We have not crafted an algorithm that can construct these linkages in a manner that respects no-overlap constraints.

\section{Conclusions and Future Work}

We have established principles that allow tight alternating torus links to be constructed, based on optimizing a multihelical unit that can be concatenated and linked. We have shown how the limiting ropelength per crossing depends on the number of concentric shells in the multihelix, and derived efficient arrangements for very large systems. When closed into torus links, our constructions have a ropelength-per-crossing that can be made smaller than any previously studied knot. We have established upper bounds that scale with the smallest possible crossing exponent of 3/4, with prefactors that approach those found from numerical shrinking.

The truly minimal multihelix may be found for larger $Q$ than 39 if an exact expression for $\tilde{H}(R)$ were available, as it would allow better a better search strategy for the minimal-length configuration than exhaustive combinatorics. Going beyond the standard helical construction, several rigid features of our concentric helical constructs and tori may be relaxed to improve their ropelength. Allowing the central rod to twist slightly would lead to a small reduction as shown by Kim et al. \cite{kim2024efficiency} but would not change the asymptotic scaling. The helices in each shell of our multihelices are arranged in circles, but a circle does not minimize the convex hull around $N$ unit disks until $N$ is in the hundreds \cite{kallrath2019packing}, and the length of the outer helices may be reduced if the layout of the cross-section is varied. 

The next improvement would likely come from relaxing the requirement that the helices be described only by a single sine and cosine. Examining tight annealings of torus knots, it is observed that the curves pack more tightly in the interior, where they are most constrained, by aligning at a steeper angle with respect to the surface of the torus, and curving into a shallower angle in the exterior where they are less constrained. Consistent with this, the planar radial coordinate of helices in tight torus links appears cycloidal rather than sinusoidal when plotted. The equation of each helix may be modified to include higher Fourier modes, to have broader maxima and minima, or to approximate a cycloid. Experimenting with these modifications yielded configurations that were a few percent tighter. This was not deemed worth pursuing in this work, but may be a subject of future work. Another potential improvement would involve taking two parallel helices and twisting them about each other in a superhelix, similar to the tightest annealings of alternating torus knots \cite{pieranski1998search}. It is unknown whether the best form of such a superhelix would be a symmetric double helix or something closer to a rod of Asclepius. Both these improvements rely on the fact that the negative Gaussian curvature at the interior of either the torus or the double helix allows closer packing of curves, suggesting a potential connection with differential geometry.

\section{Acknowledgements}
ARK is supported by the National Science Foundation, grant number 2336744.

\bibliographystyle{unsrt}
\bibliography{knotrefs}

\begin{thebibliography}{10}
\bibitem{wiki}
The caduceus image is taken from the wikipedia page on the topic and is by user rama. the rod of asclepius image is adapted from an image by user verdy p that appears on the wikipedia page.
\end{thebibliography}

\section{Appendix}
\subsection{Fit Function for $\tilde{H},\  \tilde{R}$}
For a given $N$, the radius of a cylindrical $N$-helix cannot be smaller than $1/\sin{\frac{\pi}{N}}\approx N/\pi$, below which the helices would overlap in the plane. Likewise, the height cannot be below $2N$ or the helices would overlap vertically. We have observed that when numerical data for $R$ and $H$ have their minimum values subtracted, and then are divided by $N$, all data follow a common trend independent of $N$ (Fig. 6). This common trend passes through the ideal helix values at (0.132,0.828) and follows two different power laws as $R$ and $H$ approach their minimal values. To describe the common trend, we transform the height and radius relative to their minimum and number:

\begin{equation}
    X(R,N)=\frac{R}{N}-\frac{1}{N\sin{\frac{\pi}{N}}}\approx\frac{R}{N}-\frac{1}{\pi}, \ \ \ Y(H,N)=\frac{H}{N}-2 
    \label{eq:xy}
\end{equation}

We used a four-parameter fit to the logarithm of the measured $Y$ data against $X$, and vice versa:

\begin{equation}
    Y(X)=\frac{0.0643X^{-1.8612}}{(0.1932X^{-1.8612}+1)^{0.5556}},
    \label{eq:fit}
\end{equation}
\begin{equation}
    X(Y)=\frac{0.1301Y^{-1.3239}}{(0.4165Y^{-1.3239}+1)^{0.5896}}.
\end{equation}
These can he transformed into $\tilde{R}$ and $\tilde{H}$ by inverting the definitions in Eq. \ref{eq:xy}. These fits are not necessarily upper bounds on radius and height, and can lead to slightly overlapping conditions for small heights and radii. The coefficient in the numerator may be increased slightly to ensure non-violation of the constraints. Because these fit functions are conjectural and have limited accuracy, any configuration generated from them or minimizations involving them should be checked for overlaps. These fit functions are less effective when the radius is much smaller than the ideal value, or the height is much larger. While these functions were useful in developing a better understanding of multihelix optimization, most results in this paper do not explicitly rely on them. The points in Fig. 4a were found using these functions, but they only differ by 1\% from similar values found using ``ideal'' radius and height.

\begin{figure}[ht]
    \centering
    \includegraphics[width=0.75\linewidth]{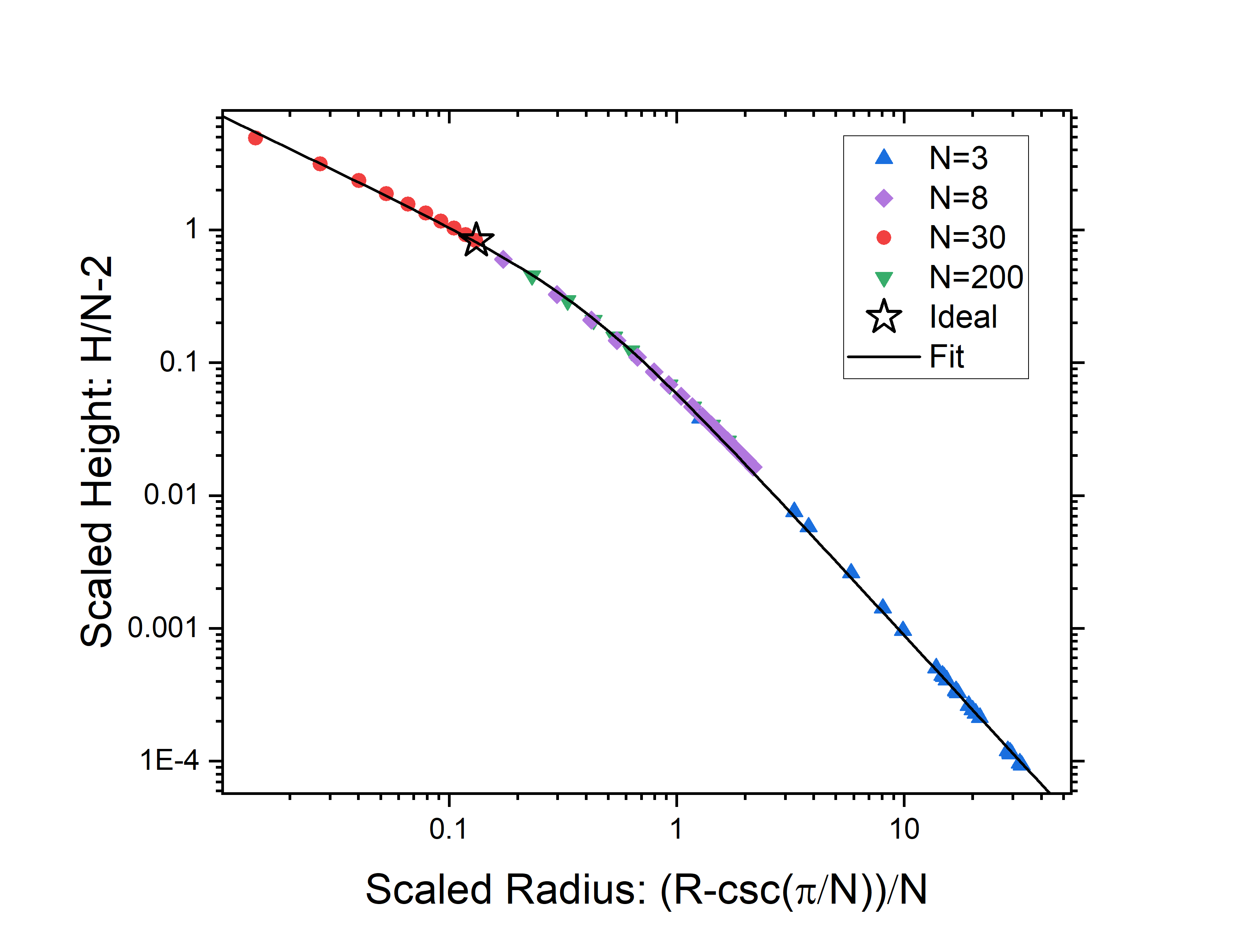}
    \caption{Scaled relationship between the height and radius of N-helices, measured numerically for different $N$. Upon rescaling, they fall upon a common curve that contains the ideal value. A four-parameter fit function, Eq. \ref{eq:fit} is overlaid. Not all data used for the fit are shown.}
    \label{fig:enter-label}
\end{figure}

\subsection{Complete table of ideal helices}
Table II in the text showed the ropelengths and ideal arrangements for a few helical constructions. Here in Table III we show the full table up to 39 helices and also include the best ratio of ropelength to the three-quarter power of crossing number for $T(3Q,Q)$ links. 

\begin{table}[ht]
\begin{tabular}{|l|l|l|l|l|}
\hline
\begin{tabular}[c]{@{}l@{}}Number of\\ Components\end{tabular} & \begin{tabular}[c]{@{}l@{}}Ropelength\\ per Twist\end{tabular} & \begin{tabular}[c]{@{}l@{}}Ropelength\\ per Crossing\end{tabular} & \begin{tabular}[c]{@{}l@{}}Best C\textasciicircum{}\{3/4\} \\ Prefactor\end{tabular} & \begin{tabular}[c]{@{}l@{}}Ideal\\ Arrangement\end{tabular} \\ \hline
2 & 14.72454 & 7.362268 & 13.8768 & {[}1,1{]} \\ \hline
3 & 30.73456 & 5.122427 & 13.4973 & {[}1,2{]} \\ \hline
4 & 49.78447 & 4.148706 & 12.5865 & {[}1,3{]} \\ \hline
5 & 75.83379 & 3.791689 & 12.4741 & {[}1,3,1{]} \\ \hline
6 & 101.8831 & 3.396104 & 13.7952 & {[}1,3,2{]} \\ \hline
7 & 127.9324 & 3.04601 & 13.0972 & {[}1,3,3{]} \\ \hline
8 & 157.7684 & 2.817293 & 11.6443 & {[}1,4,3{]} \\ \hline
9 & 185.0106 & 2.569592 & 12.6218 & {[}1,4,4{]} \\ \hline
10 & 214.4321 & 2.382579 & 12.2666 & \textit{{[}1,4,5{]}} \\ \hline
11 & 253.6763 & 2.306148 & 12.22 & \textit{{[}1,4,5,1{]}} \\ \hline
12 & 292.9205 & 2.219095 & 12.1181 & \textit{{[}1,4,5,2{]}} \\ \hline
13 & 332.1647 & 2.129261 & 11.9766 & \textit{{[}1,4,5,3{]}} \\ \hline
14 & 371.4089 & 2.040708 & 12.8455 & \textit{{[}1,4,5,4{]}} \\ \hline
15 & 410.6531 & 1.955491 & 12.613 & \textit{{[}1,4,5,5{]}} \\ \hline
16 & 459.9192 & 1.91633 & 12.04 & \textit{{[}1,5,7,3{]}} \\ \hline
17 & 500.1824 & 1.838906 & 12.5978 & \textit{{[}1,4,6,6{]}} \\ \hline
18 & 542.8946 & 1.774165 & 11.7 & \textit{{[}1,5,7,5{]}} \\ \hline
19 & 584.3822 & 1.70872 & 11.5455 & \textit{{[}1,5,7,6{]}} \\ \hline
20 & 625.8699 & 1.647026 & 11.91 & \textit{{[}1,5,7,7{]}} \\ \hline
21 & 674.2664 & 1.605396 & 12.2183 & \textbf{{[}1,5,7,8{]}} \\ \hline
22 & 727.5485 & 1.57478 & 12.1617 & \textbf{{[}1,5,7,8,1{]}} \\ \hline
23 & 780.8305 & 1.543143 & 12.0935 & \textbf{{[}1,5,7,8,2{]}} \\ \hline
24 & 834.1126 & 1.511074 & 12.0181 & \textbf{{[}1,5,7,8,3{]}} \\ \hline
25 & 887.3946 & 1.478991 & 11.9359 & \textbf{{[}1,5,7,8,4{]}} \\ \hline
26 & 940.6767 & 1.447195 & 11.8492 & \textbf{{[}1,5,7,8,5{]}} \\ \hline
27 & 993.9587 & 1.415896 & 11.7592 & \textbf{{[}1,5,7,8,6{]}} \\ \hline
28 & 1047.241 & 1.385239 & 12.4236 & \textbf{{[}1,5,7,8,7{]}} \\ \hline
29 & 1100.523 & 1.355324 & 12.3093 & \textbf{{[}1,5,7,8,8{]}} \\ \hline
30 & 1165.793 & 1.339992 & 12.2867 & \textbf{{[}1,5,7,8,8,1{]}} \\ \hline
31 & 1226.158 & 1.31845 & 12.412 & \textit{{[}1,5,8,9,8{]}} \\ \hline
32 & 1280.561 & 1.290888 & 12.298 & \textit{{[}1,5,8,9,9{]}} \\ \hline
33 & 1346.749 & 1.27533 & 12.2626 & \textit{{[}1,5,8,9,9,1{]}} \\ \hline
34 & 1412.937 & 1.259302 & 12.2219 & \textit{{[}1,5,8,9,9,2{]}} \\ \hline
35 & 1479.125 & 1.242962 & 12.1769 & \textit{{[}1,5,8,9,9,3{]}} \\ \hline
36 & 1545.313 & 1.226439 & 12.1283 & \textit{{[}1,5,8,9,9,4{]}} \\ \hline
37 & 1611.501 & 1.209836 & 12.0767 & \textit{{[}1,5,8,9,9,5{]}} \\ \hline
38 & 1677.689 & 1.193236 & 12.0227 & \textit{{[}1,5,8,9,9,6{]}} \\ \hline
39 & 1743.877 & 1.176705 & 11.97 & \textit{{[}1,5,8,9,9,7{]}} \\ \hline
\end{tabular}
\caption{Table of ropelength (per twist and crossing) of optimal concentric helices, ratio of the ropelength of a $T(3Q,Q)$ torus link to the three-quarter power of crossing number, and the ideal arrangement of helices. Italic arrangements have the height determined by the second shell (third entry), bold arrangements by the third shell (fourth entry), and normal arrangements by the first shell (second entry). }
\end{table}

\subsection{Sum Convergence}
To evaluate the sum in Eq. \ref{eq:alpha}, we can replace it with an integral over $di$ assuming T is large:
\begin{equation}
    \sqrt{8}k^2\sum_{i=1}^{T-1}i\sqrt{T^{2}+i^{2}}\approx\sqrt{8}k^2\int_{1}^{T-1}i\sqrt{T^{2}+i^{2}}di=\sqrt{8}k^2\left.\frac{1}{3}(T^2+i^2)^{3/2}\right|^{T-1}_{1}
\end{equation}
The evaluation of the integral may be simplified when T is large:
\begin{equation}
    \frac{\sqrt{8}k^2}{3}\left(\left(T^2+(T-1)^2\right)^{3/2}-\left(T^2+1^2\right)^{3/2}\right)\approx\frac{\sqrt{8}k^2}{3}\left(\sqrt{8}-1\right)T^3\approx1.72k^2T^3
\end{equation}
Note that the numerical factor is very close to but distinct from the square root of three.

In Eq. \ref{eq:upper} we take the sum of two limiting cases to establish an upper bound. To show this, we can compare the summands of the two approximations and of the exact expression, assuming ideal values for the height and inner radius:
\begin{equation}
    4\pi(i-1)+4N_s>_? \sqrt{8N_s^2 + (2\pi(R_s+2(i-1)))^2}=\sqrt{16N_s^2 + 16\pi^2(i-1)i}
\end{equation}
We factor out the 4 and square both sides:
\begin{equation}
    N_s^2+\pi^2(i-1)^2+2N_s\pi(i-1)>_?N_s^2 + \pi^2(i-1)i
\end{equation}
Subtracting the right from the left yields:
\begin{equation}
    (i-1)(\pi^2(i-1)+2N_s\pi)-(i-1)\pi^2 i>_?0\rightarrow 2N_s\pi-\pi^2>_?0
\end{equation}
This inequality holds if $i>1$ and $N_s>\pi/2$, establishing the cromulence of our upper bound.

The sum in Eq. \ref{eq:934} has a more complex integral evaluation:
\begin{equation}
    \sum_{i=1}^{\sqrt{\frac{2Q}{\pi}}}\sqrt{\frac{\pi Q}{2}}\sqrt{4\pi Q+\left(2\pi\left(\sqrt{\frac{Q}{\pi}}+2(i-1)\right)\right)^2}\approx\int_{1}^{\sqrt{\frac{2Q}{\pi}}}\sqrt{4\pi Q+\left(2\pi\left(\sqrt{\frac{Q}{\pi}}+2(i-1)\right)\right)^2}di
\end{equation}

In the large-Q limit, this evaluates to:
\begin{equation}
\sqrt{\frac{\pi}{8}}\left(2\sqrt{10+4\sqrt{2}}-\sqrt{2}+\log\frac{\sqrt{2}-1}{\sqrt{10+4\sqrt{2}}+\sqrt{8}-1} \right)Q^{3/2}\approx9.34Q^{3/2}    
\end{equation}
 The integral has a non-aesthetic evaluation that was determined by \textit{Mathematica}. Instead of showing the gigantic equation that adds little insight, we have included the code necessary to replicate it.
\begin{verbatim}
    
t=Sqrt[2*Q/Pi];n=Sqrt[Q*Pi/2];
term=Sqrt[4*Pi*Q+( (2*Pi)*(Sqrt[2]/Pi*nv+2*(i-1)))^2];
integral=Integrate[n*term,i];
evaluation=(integral/.i->t)-(integral/.i->1);
value=Limit[evaluation/Q^(3/2),Q->Infinity]
N[value]
\end{verbatim}

\subsection{Lower bound of $T(3Q,Q)$ links}
Imagine one of the components in a torus link such as that pictured in Fig. 7 is tightened, while the remaining $Q-1$ helices are not. Since each component has a linking number of 3 with the others, at its tightest it will confine $3(Q-1)$ unit disks when viewed in cross-section. Following a theorem from Cantarella et al. \cite{cantarella2002minimum}, the circumference of this link cannot be shorter than the minimum convex hull around the same number of disks. When $Q$ is large, the arrangement of disks with the smallest circumference is that of hexagonal close packing bounded by a circular hull. The radius of the circle is the related to the number of disks by their area:

\begin{equation}
    \sigma\pi R^2=3(Q-1)\pi\rightarrow R\approx\sqrt{\frac{3Q}{\sigma}}
\end{equation}
Here, $\sigma=\pi/\sqrt{12}$ is the hexagonal packing fraction of circles, and the factor of $\pi$ on the right is the area of a unit disk. The circumference of each component is $2\pi$ times the radius. Not each of the $Q$ components can be minimal in length, thus assuming that each one is minimal leads to a lower bound:
\begin{equation}
    L_{lower}=Q\times 2\pi\sqrt{\frac{3Q}{\sigma}}=2\pi\sqrt{\frac{3}{\sigma}}Q^{3/2}
\end{equation}
The crossing number is approximately $3Q^2$, and the ratio $L/C^{3/4}$ follows:
\begin{equation}
    \frac{L_{lower}}{C^{3/4}}=2\pi\sqrt{\frac{3}{\sigma}}\frac{Q^{3/2}}{(3Q^2)^{3/4}}=2\sqrt{\pi}\frac{12^{1/4}}{3^{1/4}}=\sqrt{8\pi}\approx 5.01
\end{equation}
The same logic may be applied to other $Q$; the lower bound is slightly higher for $Q=2$ (5.55), but as $Q$ gets into the hundreds the lower bound may fall below stronger pre-existing bounds. For packings below 100 disks (corresponding to $Q$ in the 30s), the minimal convex hull is not well-approximated by a circle and approximate values of computationally-determined hulls may be used \cite{kallrath2021near}. Minimal hulls around disks are computed for unit diameter (by convention) and do not take into account the offset required for radius of the surrounding component. If $N$ disks are best-encapsulated by a hull of length $L_{H_N}$, the lower bound of the ropelength of a $T(3Q,Q)$ link is $L_{lower}=Q(2L_{H_{3(Q-1)}}+2\pi)$. Applying this method to the T(4,2) Solomon link (with non-conjectured hulls) establishes a lower bound of 33.1, 18\% below the best computed upper bound \cite{ashton2011knot} and 31\% stronger than the best previously published lower bound \cite{diao2003lower}. However, since many minimal convex hulls above $N$=6 are only conjectured to be truly minimal based on computations, these lower bounds should be treated as conjectural as well. If a better convex hull is discovered for a given number of disks, it would lower the corresponding lower bound for torus links, which would mean it was not a true lower bound.

\begin{figure}[ht]
    \centering
    \includegraphics[width=0.8\linewidth]{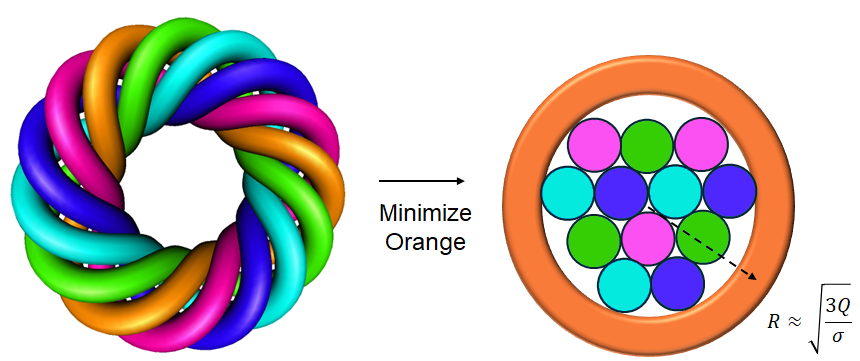}
    \caption{A T(15,5) link, and the cross section after one component has been tightened into a circle, forcing the other components into a hexagonal packing. The packed configuration is meant to be illustrative and not the true minimal packing of 12 disks.}
    \label{fig:SIhull}
\end{figure}

The hexagonally packed hull argument was used in previous work to show that Hopf linkings of two sets of untwisted, unknotted bundles of circles would have a limiting ropelength of 8.65$C^{3/4}$ \cite{klotz2021ropelength}. This was the result of both a factor-of-two error and a minor algebraic error: the prefactor should be $8\sqrt{\frac{\pi}{2}}6^{\frac{1}{4}}\approx15.7$, which is improved upon by our constructed torus links. Any reader who has made it to the end of the appendix would likely also like to know that in the original derivation used to established the current best lower bound prefactor of $(4\pi/11)^{3/4}\approx1.10$ \cite{buck1999thickness}, the 11 was an approximation that can be replaced by 10.67, increasing the lower bound to 1.13.

\end{document}